\documentclass[notitlepage,leqno,11pt]{amsart}
\usepackage{amsmath,amssymb,amsbsy,amsfonts,amsthm,latexsym,
            amsopn,amstext,amsxtra,euscript,amscd}
\usepackage{hyperref}

\begin{document}
\bibliographystyle{plain}
\newfont{\teneufm}{eufm10}
\newfont{\seveneufm}{eufm7}
\newfont{\fiveeufm}{eufm5}
%
%
\newfam\eufmfam
              \textfont\eufmfam=\teneufm \scriptfont\eufmfam=\seveneufm
              \scriptscriptfont\eufmfam=\fiveeufm
\def\bbbr{{\rm I\!R}}
\def\bbbm{{\rm I\!M}}
\def\bbbn{{\rm I\!N}}
\def\bbbf{{\rm I\!F}}
\def\bbbh{{\rm I\!H}}
\def\bbbk{{\rm I\!K}}
\def\bbbp{{\rm I\!P}}
\def\bbbone{{\mathchoice {\rm 1\mskip-4mu l} {\rm 1\mskip-4mu l}
{\rm 1\mskip-4.5mu l} {\rm 1\mskip-5mu l}}}
\def\bbbc{{\mathchoice {\setbox0=\hbox{$\displaystyle\rm C$}\hbox{\hbox
to0pt{\kern0.4\wd0\vrule height0.9\ht0\hss}\box0}}
{\setbox0=\hbox{$\textstyle\rm C$}\hbox{\hbox
to0pt{\kern0.4\wd0\vrule height0.9\ht0\hss}\box0}}
{\setbox0=\hbox{$\scriptstyle\rm C$}\hbox{\hbox
to0pt{\kern0.4\wd0\vrule height0.9\ht0\hss}\box0}}
{\setbox0=\hbox{$\scriptscriptstyle\rm C$}\hbox{\hbox
to0pt{\kern0.4\wd0\vrule height0.9\ht0\hss}\box0}}}}
\def\bbbq{{\mathchoice {\setbox0=\hbox{$\displaystyle\rm
Q$}\hbox{\raise
0.15\ht0\hbox to0pt{\kern0.4\wd0\vrule height0.8\ht0\hss}\box0}}
{\setbox0=\hbox{$\textstyle\rm Q$}\hbox{\raise
0.15\ht0\hbox to0pt{\kern0.4\wd0\vrule height0.8\ht0\hss}\box0}}
{\setbox0=\hbox{$\scriptstyle\rm Q$}\hbox{\raise
0.15\ht0\hbox to0pt{\kern0.4\wd0\vrule height0.7\ht0\hss}\box0}}
{\setbox0=\hbox{$\scriptscriptstyle\rm Q$}\hbox{\raise
0.15\ht0\hbox to0pt{\kern0.4\wd0\vrule height0.7\ht0\hss}\box0}}}}
\def\bbbt{{\mathchoice {\setbox0=\hbox{$\displaystyle\rm.
T$}\hbox{\hbox to0pt{\kern0.3\wd0\vrule height0.9\ht0\hss}\box0}}
{\setbox0=\hbox{$\textstyle\rm T$}\hbox{\hbox
to0pt{\kern0.3\wd0\vrule height0.9\ht0\hss}\box0}}
{\setbox0=\hbox{$\scriptstyle\rm T$}\hbox{\hbox
to0pt{\kern0.3\wd0\vrule height0.9\ht0\hss}\box0}}
{\setbox0=\hbox{$\scriptscriptstyle\rm T$}\hbox{\hbox
to0pt{\kern0.3\wd0\vrule height0.9\ht0\hss}\box0}}}}
\def\bbbs{{\mathchoice
{\setbox0=\hbox{$\displaystyle     \rm S$}\hbox{\raise0.5\ht0\hbox
to0pt{\kern0.35\wd0\vrule height0.45\ht0\hss}\hbox
to0pt{\kern0.55\wd0\vrule height0.5\ht0\hss}\box0}}
{\setbox0=\hbox{$\textstyle        \rm S$}\hbox{\raise0.5\ht0\hbox
to0pt{\kern0.35\wd0\vrule height0.45\ht0\hss}\hbox
to0pt{\kern0.55\wd0\vrule height0.5\ht0\hss}\box0}}
{\setbox0=\hbox{$\scriptstyle      \rm S$}\hbox{\raise0.5\ht0\hbox
to0pt{\kern0.35\wd0\vrule height0.45\ht0\hss}\raise0.05\ht0\hbox
to0pt{\kern0.5\wd0\vrule height0.45\ht0\hss}\box0}}
{\setbox0=\hbox{$\scriptscriptstyle\rm S$}\hbox{\raise0.5\ht0\hbox
to0pt{\kern0.4\wd0\vrule height0.45\ht0\hss}\raise0.05\ht0\hbox
to0pt{\kern0.55\wd0\vrule height0.45\ht0\hss}\box0}}}}
\def\bbbz{{\mathchoice {\hbox{$\sf\textstyle Z\kern-0.4em Z$}}
{\hbox{$\sf\textstyle Z\kern-0.4em Z$}}
{\hbox{$\sf\scriptstyle Z\kern-0.3em Z$}}
{\hbox{$\sf\scriptscriptstyle Z\kern-0.2em Z$}}}}
\def\ts{\thinspace}

\newtheorem{theorem}{Theorem}
\newtheorem{lemma}[theorem]{Lemma}
\newtheorem{claim}[theorem]{Claim}
\newtheorem{cor}[theorem]{Corollary}
\newtheorem{prop}[theorem]{Proposition}
\newtheorem{definition}[theorem]{Definition}
\newtheorem{remark}[theorem]{Remark}
\newtheorem{question}[theorem]{Open Question}
\newtheorem{example}[theorem]{Example}

\def\qed{\ifmmode
\squareforqed\else{\unskip\nobreak\hfil
\penalty50\hskip1em\null\nobreak\hfil\squareforqed
\parfillskip=0pt\finalhyphendemerits=0\endgraf}\fi}

\def\squareforqed{\hbox{\rlap{$\sqcap$}$\sqcup$}}

\def \C {{\mathbb C}}
\def \F {{\mathbb F}}
\def \L {{\mathbb L}}
\def \K {{\mathbb K}}
\def \Q {{\mathbb Q}}
\def \Z {{\mathbb Z}}
\def\cA{{\mathcal A}}
\def\cB{{\mathcal B}}
\def\cC{{\mathcal C}}
\def\cD{{\mathcal D}}
\def\cE{{\mathcal E}}
\def\cF{{\mathcal F}}
\def\cG{{\mathcal G}}
\def\cH{{\mathcal H}}
\def\cI{{\mathcal I}}
\def\cJ{{\mathcal J}}
\def\cK{{\mathcal K}}
\def\cL{{\mathcal L}}
\def\cM{{\mathcal M}}
\def\cN{{\mathcal N}}
\def\cO{{\mathcal O}}
\def\cP{{\mathcal P}}
\def\cQ{{\mathcal Q}}
\def\cR{{\mathcal R}}
\def\cS{{\mathcal S}}
\def\cT{{\mathcal T}}
\def\cU{{\mathcal U}}
\def\cV{{\mathcal V}}
\def\cW{{\mathcal W}}
\def\cX{{\mathcal X}}
\def\cY{{\mathcal Y}}
\def\cZ{{\mathcal Z}}
\newcommand{\rmod}[1]{\: \mbox{mod}\: #1}

\def\tcN{\cN^\mathbf{c}}
\def\F{\mathbb F}
\def\Tr{\operatorname{Tr}}
\def\mand{\qquad \mbox{and} \qquad}
\renewcommand{\vec}[1]{\mathbf{#1}}
\def\eqref#1{(\ref{#1})}
\newcommand{\ignore}[1]{}
\hyphenation{re-pub-lished}
\parskip 1.5 mm
\def\lln{{\mathrm Lnln}}
\def\Res{\mathrm{Res}\,}
\def\F{{\bbbf}}
\def\Fp{\F_p}
\def\fp{\Fp^*}
\def\Fq{\F_q}
\def\ff{\F_2}
\def\ffn{\F_{2^n}}
\def\K{{\bbbk}}
\def \Z{{\bbbz}}
\def \N{{\bbbn}}
\def\Q{{\bbbq}}
\def \R{{\bbbr}}
\def \P{{\bbbp}}
\def\Zm{\Z_m}
\def \Um{{\mathcal U}_m}
\def \Bf{\frak B}
\def\Km{\cK_\mu}
\def\va {{\mathbf a}}
\def \vb {{\mathbf b}}
\def \vc {{\mathbf c}}
\def\vx{{\mathbf x}}
\def \vr {{\mathbf r}}
\def \vv {{\mathbf v}}
\def\vu{{\mathbf u}}
\def \vw{{\mathbf w}}
\def \vz {{\mathbfz}}
\def\\{\cr}
\def\({\left(}
\def\){\right)}
\def\fl#1{\left\lfloor#1\right\rfloor}
\def\rf#1{\left\lceil#1\right\rceil}
\def\flq#1{{\left\lfloor#1\right\rfloor}_q}
\def\flp#1{{\left\lfloor#1\right\rfloor}_p}
\def\flm#1{{\left\lfloor#1\right\rfloor}_m}
\def\Al{{\sl Alice}}
\def\Bob{{\sl Bob}}
\def\Or{{\mathcal O}}
\def\inv#1{\mbox{\rm{inv}}\,#1}
\def\invM#1{\mbox{\rm{inv}}_M\,#1}
\def\invp#1{\mbox{\rm{inv}}_p\,#1}
\def\Ln#1{\mbox{\rm{Ln}}\,#1}
\def \nd {\,|\hspace{-1.2mm}/\,}
\def\ord{\mu}
\def\E{\mathbf{E}}
\def\Cl{{\mathrm {Cl}}}
\def\epp{\mbox{\bf{e}}_{p-1}}
\def\ep{\mbox{\bf{e}}_p}
\def\eq{\mbox{\bf{e}}_q}
\def\bm{\bf{m}}
\newcommand{\floor}[1]{\lfloor {#1} \rfloor}
\newcommand{\comm}[1]{\marginpar{
\vskip-\baselineskip
\raggedright\footnotesize
\itshape\hrule\smallskip#1\par\smallskip\hrule}}
\def\rem{{\mathrm{\,rem\,}}}
\def\dist {{\mathrm{\,dist\,}}}
\def\etal{{\it et al.}}
\def\ie{{\it i.e. }}
\def\veps{{\varepsilon}}
\def\eps{{\eta}}
\def\ind#1{{\mathrm {ind}}\,#1}
               \def \MSB{{\mathrm{MSB}}}
\newcommand{\abs}[1]{\left| #1 \right|}

\title{ Lucas Numbers with the Lehmer property }
%
\author{
{\sc Bernadette~Faye}\and {\sc Florian~Luca}
}
\address{
Ecole Doctorale de Mathematiques et d'Informatique \newline
Universit\'e Cheikh Anta Diop de Dakar \newline
BP 5005, Dakar Fann, Senegal and\newline
School of Mathematics, University of the Witwatersrand \newline
Private Bag X3, Wits 2050, Johannesburg, South Africa
}
\email{bernadette@aims-senegal.org}

\address{
School of Mathematics, University of the Witwatersrand \newline
Private Bag X3, Wits 2050, Johannesburg, South Africa \newline
}
\email{Florian.Luca@wits.ac.za}

\maketitle

\begin{abstract} A composite positive integer $n$ is \textit{Lehmer}  if $\phi(n)$ divides $n-1$, where $\phi(n)$ is the Euler's totient function. No Lehmer number is known, nor has it been proved that they don't exist. In 2007, the second author \cite{FL1} proved that there is no Lehmer number
 in the Fibonacci sequence. In this paper, we adapt the method from \cite{FL1} to show that there is no Lehmer number in the companion Lucas sequence of the Fibonacci sequence $(L_n)_{n\ge 0}$ given by $L_0=2, L_1=1$ and $L_{n+2}=L_{n+1}+L_n$ for all $n\geq0$.
\end{abstract}

\section{Introduction}
Let $\phi(n)$ be the Euler function of a positive  integer $n$.
Recall that if $n$ has the prime factorization
$$
n=p_1^{\alpha_1} p_2^{\alpha_2}\cdots p_k^{\alpha_k},
$$
then
$$
\phi(n)=(p_1-1)p_1^{\alpha_1-1} (p_2-1)p_2^{\alpha_2-1}\cdots (p_k-1)p_k^{\alpha_k-1}.
$$
Lehmer \cite{Leh} conjectured that if $\phi(n)\mid n-1$ then $n$ is a prime. To this day, the conjecture remains open. Counterexamples to Lehmer's conjecture have been dubbed {\it Lehmer numbers}.
Several people worked on getting larger and larger lower bounds on a potential Lehmer number. For a positive integer $m$, we write $\omega(m)$ for the number of  distinct prime factors of $m$. Lehmer himself proved that if $N$ is Lehmer, then $\omega(N)\geq 7$. This has been improved by Cohen and Hagis \cite{coh} to $\omega(N)\geq 14.$ The current record $\omega(N)\ge 15$ is due to Renze \cite{joh}. If additionally $3\mid N$, then $\omega(N)\geq 40\cdot 10^{6}$ and $N>10^{36\cdot 10^{7}}.$

Not succeeding in proving that there are no Lehmer numbers, some researchers have settled for the more modest goal of proving that
there are no Lehmer numbers in certain interesting subsequences of positive integers. For example, in \cite{FL1}, Luca proved that there is no Fibonacci number which is Lehmer. In \cite{GL}, it is shown that there is no Lehmer number in the sequence of Cullen numbers $\{C_n\}_{n\ge 1}$ of general term $C_n=n2^n+1$, while in \cite{Dajune} the same conclusion is shown to hold for generalized Cullen numbers. In \cite{L1},
it is shown that there is no Lehmer number of the form $(g^n-1)/(g-1)$ for any $n\ge 1$ and integer $g \in [2,1000]$.

\medskip

Here, we apply the same argument as in \cite{FL1}, to the Lucas sequence companion of the Fibonacci sequence given by $L_0=2, L_1=1$ and $L_{n+2}=L_{n+1}+L_n$ for all $n\geq0$. Putting $(\alpha,\beta)=((1+{\sqrt{5}})/2,(1-{\sqrt{5}})/2)$
for the two roots of the characteristic equation $x^2-x-1=0$ of the Lucas sequence, the Binet formula
\begin{equation}
\label{eq:1}
L_n=\alpha^n +\beta^n\qquad {\text{\rm holds~for~all}}\qquad n\ge 0.
\end{equation}
There are several relations among Fibonacci and Lucas numbers which are well-known and can be proved using the Binet formula \eqref{eq:1} for the Lucas numbers and its analog
$$
F_n=\frac{\alpha^n-\beta^n}{\alpha-\beta}\qquad {\text{\rm for~all}}\qquad n\ge 0
$$
for the Fibonacci numbers. Some of them which are useful for us are
\begin{equation}
\label{r1}
L_n^{2}-5F_n^{2}=4(-1)^n,
\end{equation}
\begin{equation}
\label{eq:2}
L_n=L_{n/2}^2-2(-1)^{n/2}\quad {\text{\rm valid~for~ all~ even}}\quad n,
\end{equation}
whereas for odd $n$
\begin{equation}
\label{eq:3}
L_n-1=\left\{ \begin{matrix}
5F_{(n+1)/2} F_{(n-1)/2} & {\text{\rm if}} & n\equiv 1\pmod 4;\\
L_{(n+1)/2} L_{(n-1)/2} & {\text{\rm if}} & n\equiv 3\pmod 4.
\end{matrix}\right.
\end{equation}
Our result is the following:
\begin{theorem}
There is no Lehmer number in the Lucas sequence.
\end{theorem}

\section{Proof}
Assume that $L_n$ is Lehmer for some $n$.
Clearly, $L_n$ is odd and $\omega(L_n)\ge 15$ by the main result from \cite{joh}. The product of the first $15$ odd primes exceeds $1.6\times 10^{19}$, so $n\ge 92$. Furthermore,

\begin{equation}
\label{eq:15}
2^{15}\mid 2^{\omega(L_n)}\mid \phi(L_n)\mid L_n-1.
\end{equation}
\medskip

If $n$ is even, formula \eqref{eq:2} shows that $L_n-1=L_{n/2}^2+1$ or $L_{n/2}^2-3$ and numbers of the form $m^2+1$ or $m^2-3$ for some integer $m$ are never multiples of $4$, so divisibility \eqref{eq:15} is impossible. If $n\equiv 3\pmod 8$, relations \eqref{eq:3} and \eqref{eq:15} show that $2^{15}\mid L_{(n+1)/2}L_{(n-1)/2}$. This is also impossible since no member of the Lucas sequence is a multiple of $8$, fact which can be easily proved by listing its first $14$ members modulo $8$:
$$
2,~1,~3,~4,~7,~3,~2,~5,~7,~4,~3,~7,~2,~1,
$$
and noting that we have already covered the full period of $\{L_m\}_{m\ge 0}$ modulo $8$ (of length $12$) without having reached any zero.

\medskip

So, we are left with the case when $n\equiv 1\pmod{4}.$

Let us write
$$
n=p_1^{\alpha_1}\cdots p_k^{\alpha_k},
$$
with $p_1<\cdots<p_k$ odd primes and $\alpha_1,\ldots,\alpha_k$ positive integers. If $p_1=3$, then $L_n$ is even, which is not the case. Thus, $p_1\geq 5$.

\medskip

Here, we use the argument from \cite{FL1}  to bound $p_1$. Since most of the details are similar, we only sketch the argument. Let $p$ be any prime factor of $L_n$. Reducing formula \eqref{eq:1} modulo $p$ we get that $-5F_n^2\equiv -4\pmod p$. In particular, $5$ is a quadratic residue modulo $p$, so by Quadratic Reciprocity also $p$ is a quadratic residue modulo $5$. Now let $d$ be any divisor of $n$ which is a multiple of $p_1$.
By Carmichael's Primitive Divisor Theorem for the Lucas numbers (see \cite{RD}), there exists a primitive prime $p_d\mid L_d$, such that $p_d\nmid L_{d_1}$ for all positive $d_1<d.$ Since $n$ is odd and $d\mid n$, we have $L_d\mid L_n$, therefore $p_d\mid L_n$. Since $p_d$ is primitive for $L_d$ and a quadratic residue modulo $5$, we have $p_d\equiv 1\pmod d$ (if $p$ were not a quadratic residue modulo $5$, then we would have had that $p_d\equiv -1\pmod 5$, which is less useful for our problem). In particular,
\begin{equation}
\label{eq:p1}
p_1\mid d\mid p_d-1\mid \phi(L_n).
\end{equation}
Collecting the above divisibilities \eqref{eq:p1} over all divisors $d$ of $n$ which are multiples of $p_1$ and using \eqref{eq:3}, we have
\begin{equation}
\label{eq:5}
p_1^{\tau(n/p_1)} \mid \phi(L_n)\mid L_n-1\mid 5 F_{(n-1)/2}F_{(n+1)/2}.
\end{equation}
In the above, $\tau(m)$ is the number of divisors of $m$. If $p_1=5$, then
$5\mid n$, therefore $5\nmid F_{(n\pm 1)/2}$ because a Fibonacci number $F_m$ is a multiple of $5$ if and only if its index $m$ is a multiple of $5$. Thus, $\tau(n/p_1)=1$, so $n=p_1$, which is impossible since $n>92$.

Assume now that $p_1>5$. Since
$$
\gcd(F_{(n+1)/2}, F_{(n-1)/2})=F_{\gcd((n+1)/2,(n-1)/2)}=F_1=1,
$$
divisibility relation \eqref{eq:5} shows that $p_1^{\tau(n/p_1)}$
divides $F_{(n+\varepsilon)/2}$ for some $\varepsilon\in \{\pm 1\}$. Let $z(p_1)$ be the order of appearance of $p_1$ in the Fibonacci sequence, which is the minimal positive integer $\ell$ such that $p_1\mid F_{\ell}$. Write
\begin{equation}
\label{eq:Wall}
F_{z(p_1)}=p_1^{e_{p_1}} m_{p_1},
\end{equation}
where $m_{p_1}$ is coprime to $p_1$. It is known that $p_1\mid F_k$ if and only if $z(p_1)\mid k$. Furthermore, if $p_1^t\mid F_k$ for some $t>e_{p_1}$, then necessarily $p_1\mid k$. Since for us $(n+\varepsilon)/2$ is not a multiple of $p_1$ (because $n$ is a multiple of $p_1$), we get that $\tau(n/p_1)\le e_{p_1}$. In particular, if $p_1=7$, then $e_{p_1}=1$, 
so $n=p_1$, which is false since $n>92$. So, $p_1\ge 11$. We now follow along the argument from \cite{FL1} to get that
\begin{equation}
\label{eq:30}
\tau(n)\leq 2\tau(n/p_1)\leq \frac{(p_1+1)\log\alpha}{\log p_1}.
\end{equation}
Further, since $(L_n-1)/\phi(L_n)$ is an integer larger than $1$, we have
\begin{equation}
\label{r:1}
2<\frac{L_n}{\phi(L_n)}\leq \prod_{p\mid L_n}\(1+\frac{1}{p-1}\)<\exp\(\sum_{p\mid L_n}\frac{1}{p-1}\),
\end{equation}
or
\begin{equation}
\label{eq:4}
\log 2 \leq \sum_{p \mid L_n}\frac{1}{p-1}.
\end{equation}
Letting for a divisor $d$ of $n$ the notation ${\mathcal P_d}$ stand for the set of primitive prime factors of $L_d$, the argument from \cite{FL1} gives
\begin{equation}
\label{eq:imp}
\sum_{p \in \mathcal{P}_d}\frac{1}{p-1} \leq \frac{0.9}{d} + \frac{2.2\log\log d}{d}.
\end{equation}
Since the function $x\mapsto (\log\log x)/x$ is decreasing for $x>10$ and all divisors $d>1$ of $n$ satisfy $d>10$, we have, using \eqref{eq:30}, that
\begin{eqnarray}
\label{eq:8}
\sum_{p\mid L_n}\frac{1}{p-1}&=& \sum_{d\mid n}\sum_{p\in \mathcal{P}_d}\frac{1}{p-1} \leq \sum_{\substack{d\mid n\\ d>1}}\(\frac{0.9}{d} + \frac{2.2\log\log d}{d}\)\\
&\leq & \(\frac{0.9}{p_1} + \frac{2.2\log\log p_1}{p_1}\)\tau(n) \nonumber \\
&\leq & (\log\alpha)\frac{(p_1+1)}{\log p_1}\cdot \(\frac{0.9}{p_1} + \frac{2.2\log\log p_1}{p_1}\),\nonumber
\end{eqnarray}
which together with inequality \eqref{eq:4} leads to
\begin{equation}
\label{eq:9}
\log p_1\leq \frac{(\log\alpha)}{\log 2}\left(\frac{p_1+1}{p_1}\right)(0.9 + 2.2\log\log p_1).
\end{equation}
The above inequality \eqref{eq:9} implies $p_1<1800$.
Since $p_1<10^{14}$, a calculation of McIntosh and Roettger \cite{MC} shows that $e_{p_1}=1$. Thus, $\tau(n/p_1)=1$, therefore $n=p_1.$ Since $n\ge 92$, we have $p_1\ge 97$. Going back to the inequalities \eqref{eq:4} and \eqref{eq:imp}, we get
$$
\log 2<\frac{0.9}{p_1}+\frac{2.2 \log\log p_1}{p_1},
$$
which is false for $p_1\ge 97$. The theorem is proved.

\section*{Acknowledgments} B. F. thanks OWSD and Sida (Swedish International Development Cooperation Agency) for a scholarship during her Ph.D. studies at Wits.


\end{document}